\documentclass[psamsfonts]{amsart}
\usepackage{amsfonts}
\usepackage{amsmath,amsthm}

 \def\CD{{\mathcal D}}
 
 \def\CH{{\mathcal H}}

 \def\CP{{\mathcal P}}

 \def\CV{{\mathcal V}}
 \def\CW{{\mathcal W}}

 \def\NN{{\mathbb N}}

 \def\RR{{\mathbb R}}
 \def\ZZ{{\mathbb Z}}
        \def\diag{\operatorname{diag}}
        \def\proj{\operatorname{proj}}

\newcommand{\wh}{\widehat}

\newtheorem{theorem}{ \noindent {Theorem}}[section]
\newtheorem{corollary}[theorem]{\noindent {Corollary}}
\newtheorem{proposition}[theorem]{ \noindent {Proposition}}
\newtheorem{definition}[theorem]{\noindent {Definition}}

\begin{document}

\title[generalized translation operator in several variables]
{Generalized translation operator and approximation in several variables}
\author{ Yuan Xu}
\address{Department of Mathematics\\ University of Oregon\\
    Eugene, Oregon 97403-1222.}\email{yuan@math.uoregon.edu}

\date{September 10, 2003}
\thanks{Work supported in part by the National Science Foundation 
under Grant DMS-0201669}
                                                    
\begin{abstract}
Generalized translation operators for orthogonal expansions with respect
to families of weight functions on the unit ball and on the standard simplex 
are studied. They are used to define convolution structures and modulus of 
smoothness for these regions, which are in turn used to characterize the 
best approximation by polynomials in the weighted $L^p$ spaces. In one 
variable, this becomes the generalized translation operator for the 
Gegenbauer polynomial expansions.
\end{abstract}

\maketitle

\section{Introduction} 
\setcounter{equation}{0}

Let $w_\lambda$ denote the weight function $w_\lambda(t)=(1-t^2)^{\lambda-1/2}$
on $[-1,1]$. Let $b_\lambda$ be the normalization constant of $w_\lambda$,
$b_\lambda^{-1} = \int_{-1}^1 w_\lambda(s)ds$. The orthogonal polynomials with
respect to $w_\lambda$ are the Gegenbauer polynomials $C_n^{\lambda}$. 
The generalized translation operator with respect to $w_\lambda$ is defined by
\begin{equation}
T_s f(t) = b_{\lambda-1/2} \int_{-1}^1 f\left(s t + 
  u \sqrt{1-s^2}\sqrt{1-t^2}\right)
     (1-u^2)^{\lambda-1} du. \label{eq:1.1}
\end{equation}
It plays the role of translation for the trigonometric series and can be 
used to define a convolution structure $f\star g$ for $f, g \in 
L^1(w_\lambda,[-1,1])$,
\begin{equation}
  (f\star g)(t)  = b_\lambda \int_{-1}^1 f(s) T_t g(s) w_\lambda(s)ds, 
\label{eq:1.2}
\end{equation}
as introduced by Gelfand \cite{Ge} and Bochner \cite{Bo}. The convolution and
the generalized translation operator have been used to study Fourier orthogonal
expansions in Gegenbauer polynomials (see, for example, 
\cite{AW,Bav,BBP,Butz,P,Raf,SW}).
Using the product formula of the Gegenbauer polynomials, the generalized 
translation operator can also be defined by the equation
\begin{equation}
 T_s C_n^\lambda(t) = \frac{C_n^\lambda(s)}{C_n^\lambda(1)} C_n^\lambda(t), 
\qquad n \ge 0. \label{eq:1.3}
\end{equation}

The purpose of this paper is to study the generalized translation operator for
weight functions defined on the unit ball $B^d = \{x: \|x\|\le 1\}\subset
\RR^d$ and on the standard simplex 
$$
T^d = \{x \in \RR^d:x_1 \ge 0, \ldots, x_d \ge 0, 1-|x| \ge 0\}, \qquad
 |x| = x_1+\ldots + x_d,
$$
in $\RR^d$ and use them to study weighted approximation and orthogonal 
expansions in several variables. 
To define the weight functions, we start from the reflection invariant weight 
function considered by Dunkl \cite{D1}. 

For a nonzero vector $v \in \RR^d$, 
let $\sigma_v$ denote the reflection with respect to the hyperplane 
perpendicular to $v$; that is, $x \sigma_v : = x - 2 (\langle x,v \rangle /
\|v\|^2) v$, $x \in \RR^d$, where $\langle x,y\rangle$ denote the usual 
Euclidean inner product and $\|x\|$ denote the usual Euclidean norm $\|x\|^2 
= \langle x, x\rangle$. The weight function $h_\kappa$ is defined by 
\begin{equation}
h_\kappa(x) = \prod_{v \in R_+} |\langle x, v\rangle|^{\kappa_v}, \qquad 
   x \in \RR^d, \label{eq:1.4}
\end{equation}
in which $R_+$ is a fixed positive root system of $\RR^d$, normalized so 
that $\langle v, v \rangle =2$ for all $v \in R_+$, and $\kappa$ is a 
nonnegative multiplicity function $v \mapsto \kappa_v$ defined on $R_+$ with 
the property that $\kappa_u = \kappa_v$ whenever $\sigma_u$ is conjugate to 
$\sigma_v$ in the reflection group $G$ generated by the reflections 
$\{\sigma_v:v \in R_+\}$. Then $h_\kappa$ is invariant under the reflection 
group $G$, a subgroup of the orthogonal group. The simplest example 
is given by the case $G=\ZZ_2^d$ for which $h_\kappa$ is just the product 
weight function
\begin{equation}
  h_\kappa (x)  = \prod_{i=1}^d |x_i|^{\kappa_i}, \qquad \kappa_i \ge 0. 
\label{eq:1.5}
\end{equation}    
Other examples include weight functions invariant under the symmetric group
and the hyperoctahedral group,
$$
\prod_{1 \le i<j \le d} |x_i- x_j|^{\kappa} \qquad \hbox{and} \qquad
\prod_{i=1}^d |x_i|^{\kappa_0}\prod_{1 \le i<j \le d}|x_i^2- x_j^2|^{\kappa},
$$
respectively. 

The weight functions on the unit ball $B^d$ that we shall consider are of 
the form
\begin{equation}
 W_{\kappa,\mu}^B(x) = h_\kappa^2(x) (1-\|x\|^2)^{\mu-1/2}, \qquad x \in B^d,
\label{eq:1.6}
\end{equation}
where $\mu \ge 0$ and $h_\kappa$ is a reflection invariant weight function 
as in (\ref{eq:1.4}), and the weight functions on the simplex that we shall 
consider are of the form
\begin{equation}
W_{\kappa,\mu}^T(x) = h_\kappa^2(\sqrt{x_1}, \ldots,\sqrt{x_d})
   (1-|x|)^{\mu-1/2} /\sqrt{x_1 \cdots x_d},  
\label{eq:1.7}
\end{equation}
where $\mu \ge 0$ and $h_\kappa$ is a reflection invariant weight function 
as in (\ref{eq:1.4}), and we assume that $h_\kappa$ is even in each of its 
variables (for example, weight functions invariant under $\ZZ_2^d$ and the 
hyperoctahedral group on $\RR^d$). These include the classical weight 
functions on these domains, which are
\begin{equation}
W_\mu^B(x) = (1-\|x\|^2)^{\mu-1/2}, \qquad x \in B^d, \label{eq:1.8}
\end{equation}
on the unit ball (taking $h_\kappa (x) =1$) and 
\begin{equation}
W_\kappa^T(x) = x_1^{\kappa_1-1/2} \cdots x_d^{\kappa_d-1/2} 
     (1-|x|)^{\kappa_{d+1}-1/2},\qquad x \in T^d, \label{eq:1.9}
\end{equation} 
on the simplex (taking $h_\kappa(x) = \prod_{i=1}^d|x_i|^\kappa_i$ and
$\kappa_{d+1} = \mu$). For $d =1$, $W_\kappa^T(x)$ is the Jacobi
weight function on the interval $[0,1]$. 

The orthogonal structures for $W_{\kappa,\mu}^B$ on the ball and for 
$W_{\kappa,\mu}^T$ on the simplex are closely related to the orthogonal
structure of $h$-harmonics on the unit sphere $S^d = \{x: \|x\| =1\}$ of 
$\RR^{d+1}$. Our study of the generalized translation
operators relies on that of the weighted spherical means, studied in
\cite{X02a,X02b}, which are the generalizations of the usual spherical means 
\begin{equation}
  T_\theta f(x) = \frac{1}{\sigma_d (\sin \theta)^d} 
   \int_{\langle x, y\rangle = \cos \theta} f(y) d\omega(y), \label{eq:1.10}
\end{equation}
where $\sigma_d = \int_{S^d} d\omega = 2 \pi^{(d+1)/2}/\Gamma((d+1)/2)$ is 
the surface area of $S^d$. The weighted spherical means are defined implicitly
via an integral relation. In \cite{X02b} the weighted means were used to 
define a modulus of smoothness, which was shown to be equivalent to a 
K-functional and used to characterize the weighted $L^p$ best approximation 
by polynomials. The similar K-functional was also defined and used to 
characterize the weighted best approximation on $B^d$ and on $T^d$, but the 
modulus of smoothness was not defined since the analog of the spherical means
for $B^d$ and $T^d$ seemed to be artificial. It has been realized only 
recently in \cite{X03} that the analog of the spherical means for $B^d$ and
$T^d$, as generalized translation operators on these regions, are of 
interests. Since the main purpose of \cite{X03} is to define weighted maximal 
functions and use them to prove results on almost everywhere convergence, the 
generalized translation operators themselves were not studied there. We 
complete this circle of ideas in the present paper. 

One of our results gives an explicit integral formula for the generalized 
translation operator with respect to the classical weight function $W_\mu^B$
in (\ref{eq:1.8}) on the unit ball (see (\ref{eq:3.9}), which extends the 
formula (\ref{eq:1.1}) to several variables. No integral formula is known for
any other weight functions. 

The paper is organized as follows. In the next section we recall the 
background and results for $h$-harmonics and the weighted spherical means. 
The results on the unit ball are presented in Section 3 and the results on 
the simplex appear in Section 4.

\section{Weighted spherical means and weighted approximation on $S^{d-1}$}
\setcounter{equation}{0}

Let $h_\kappa$ be the reflection invariant weight function defined in 
(\ref{eq:1.4}). We denote by $a_\kappa$ the normalization constant of 
$h_\kappa$, $a_\kappa^{-1}  = \int_{S^{d-1}} h_\kappa^2(y) d\omega$, and denote
by $L^p(h_\kappa^2)$, $1\le p\le \infty$, the space of functions defined on 
$S^{d-1}$ with the finite norm
$$
\|f\|_{\kappa,p} := \Big(a_\kappa \int_{S^{d-1}} |f(y)|^p h_\kappa^2(y)
d\omega(y) \Big)^{1/p}, \qquad  1 \le p < \infty,
$$
and for $p = \infty$ we assume that $L^\infty$ is replaced by $C(S^{d})$, the 
space of continuous functions on $S^{d-1}$ with the usual uniform norm 
$\|f\|_\infty$. The case $\kappa \equiv 0$ corresponds to the usual
(unweighted)  $L^p$ space on $S^{d-1}$.

\subsection{Background}
The essential ingredient of the theory of $h$-harmonics is a family of 
first-order differential-difference operators, $\CD_i$, called 
Dunkl's operators, which generates a commutative algebra; these operators are 
defined by (\cite{D1})
$$
  \CD_i f(x) = \partial_i f(x) + \sum_{v \in R_+} k_v 
    \frac{f(x) -  f(x \sigma_v)} {\langle x, v\rangle}
        \langle v,\varepsilon_i\rangle, \qquad 1 \le i \le d, 
$$
where $\varepsilon_1, \ldots, \varepsilon_{d}$ are the standard unit vectors
of $\RR^d$. The $h$-Laplacian is defined by $\Delta_h=\CD_1^2 + \ldots + 
\CD_{d}^2$ and it plays the role similar to that of the ordinary Laplacian. 
Let $\CP_n^{d}$ denote the subspace of homogeneous polynomials of degree $n$ 
in $d$ variables. An $h$-harmonic polynomial $P$ of degree $n$ is a 
homogeneous polynomial $P \in \CP_n^{d}$ such that $\Delta_h P  =0$. 
Furthermore, let $\CH_n^{d}(h_\kappa^2)$ denote the space of $h$-harmonic 
polynomials of degree $n$ in $d$ variables and define 
$$
 \langle f, g\rangle_\kappa : = a_\kappa \int_{S^{d-1}} f(x) g(x) 
    h^2_\kappa(x) d\omega(x). 
$$ 
Then $\langle P,Q\rangle_{\kappa} = 0$ for $P \in \CH_n^{d}(h_\kappa^2)$ and 
$Q \in \Pi_{n-1}^{d}$, where $\Pi_n^{d}$ denote the space of polynomials 
of degree at most $n$ in $d$ variables. The spherical $h$-harmonics are the 
restriction of $h$-harmonics on the unit sphere. It is known that $\dim 
\CH_n(h_\kappa^2) = \dim \CP_n^d -\dim \CP_{n-2}^d$ with $\dim \CP_n^d = 
\binom{n+d -1}{d}$. 

In terms of the polar coordinates $y =ry'$, $r = \|y\|$, the $h$-Laplacian
operator $\Delta_h$ takes the form (\cite{X01b})
$$ 
\Delta_h = \frac{\partial^2}{\partial r^2} + \frac{2 \lambda_\kappa+1}{r} 
    \frac{\partial}{\partial r} + \frac{1}{r^2} \Delta_{h,0}, 
$$
where $\Delta_{h,0}$ is the (Laplace-Beltrami) operator on the sphere,
and throughout this paper, we fix the value of $\lambda_\kappa$ as 
\begin{equation}
  \lambda := \lambda_\kappa = \gamma_\kappa+ \frac{d-1}{2} \qquad 
 \hbox{with} \qquad  \gamma_\kappa =  \sum_{v\in R_+} \kappa_v.
\label{eq:2.4}
\end{equation}
Applying $\Delta_h$ to $h$-harmonics $Y \in \CH_n(h_\kappa^2)$ with $Y(y) = 
r^nY(y')$ shows that spherical $h$-harmonics are eigenfunctions of 
$\Delta_{h,0}$; that is,
\begin{equation}
 \Delta_{h,0} Y(x) = -n (n+ 2 \lambda_\kappa) Y(x), 
      \qquad x \in S^{d-1}, \quad Y \in \CH_n^d(h_\kappa^2). \label{eq:2.1}
\end{equation}
For further background materials, see \cite{D1,DX} and the references 
in \cite{DX}.

The standard Hilbert space theory shows that
$$
L^2(h_\kappa^2) = \sum_{n=0}^\infty\bigoplus \CH_n^{d}(h_\kappa^2).
$$
That is, with each $f\in L^2(h_\kappa^2)$ we can associate its $h$-harmonic
expansion
$$
  f(x) = \sum_{n=0}^\infty Y_n(h_\kappa^2;f,x), \qquad x \in S^{d-1},
$$
in $L^2(h_\kappa^2)$ norm. For the surface measure ($\kappa =0$), such a
series is called the Laplace series (cf. \cite[Chapt. 12]{Er}). The orthogonal
projection $Y_n(h_\kappa^2): L^2(h_\kappa^2) \mapsto \CH_n^{d}(h_\kappa^2)$
takes the form
\begin{equation}
 Y_n(h_\kappa^2;f,x):=
   \int_{S^{d-1}} f(y) P_n(h_\kappa^2;x,y) h_\kappa^2(y)\, d\omega(y),
\label{eq:2.2}
\end{equation}
where the kernel $P_n(h_\kappa^2;x,y)$ is the reproducing kernel of the space 
of $h$-harmonics $\CH_n^d(h_\kappa^2)$ in $L^2(h_\kappa^2)$. The kernel 
$P_n(h_\kappa^2;x,y)$ enjoys a compact formula in terms of the intertwining 
operator between the commutative algebra generated by the partial derivatives 
and the one generated by Dunkl's operators. This operator, $V_\kappa$, is 
linear and it is determined uniquely by
$$
 V_\kappa \CP_n^d \subset \CP_n^d, \qquad V_\kappa 1=1,
    \qquad \CD_i V_\kappa = V_\kappa \partial_i,  \qquad 1 \le i \le d.
$$
The formula of the reproducing kernel for $\CH_n^{d}(h_\kappa^2)$ is 
given in terms of the Gegenbauer polynomials 
\begin{equation}
P_n(h_\kappa^2;x,y) = \frac{n+\lambda_\kappa}{\lambda_\kappa}
   V_\kappa[C_n^{\lambda_\kappa} (\langle \cdot, y \rangle )](x).
\label{eq:2.3}
\end{equation}

An explicit formula of $V_\kappa$ is known only in the case of symmetric group
$S_3$ for three variables and in the case of the abelian group $\ZZ_2^{d}$. 
In the latter case, $V_\kappa$ is an integral operator,
\begin{equation} 
  V_\kappa f(x) = c_\kappa
       \int_{[-1,1]^{d}} f(x_1 t_1, \ldots,x_{d} t_{d}) 
        \prod_{i=1}^{d} (1+t_i)(1-t_i^2)^{\kappa_i -1} d t,
\label{eq:2.5}
\end{equation}
where $c_\kappa$ is the normalization constant determined by $V_\kappa 1 =1$,
$c_\kappa =   c_{\kappa_1} \ldots c_{\kappa_{d}}$ and 
$c_r^{-1} =\int_{-1}^1 (1-t^2)^{r-1}dt$.
If some $\kappa_i =0$, then the formula holds under the limit relation
$$
 \lim_{\lambda \to 0} c_\lambda \int_{-1}^1 f(t) (1-t)^{\lambda -1} dt
  = [f(1) + f(-1)] /2.
$$
One important property of the intertwining operator is that it is positive
(\cite{Ros}), that is, $V_\kappa p \ge 0$ if $p \ge 0$. One can also study
the dual of this operator, as in \cite{T}.

\subsection{Weighted spherical means, convolution and approximation}
We recall the results developed in \cite{X02a,X02b}, some of which will be
needed later and others are cited to show what can be expected in the cases
$B^d$ and $T^d$. For $f \in L^p(h_\kappa^2)$ and $g \in L^1(w_\lambda;[-1,1])$,
we define a sort of convolution 
\begin{equation} 
(f\star_\kappa g)(x): = a_\kappa \int_{S^{d-1}} f(y) 
  V_\kappa[g(\langle x,\cdot\,\rangle)](y) h_\kappa^2(y) d\omega.\label{eq:2.6}
\end{equation}
For the surface measure ($h_\kappa(x) =1$ and $V_\kappa = id$), this is called 
the spherical convolution in \cite{CZ}, and it has been used by many authors,
see, for example, \cite{Berg,BL,LW,LN,P,Rus}. It satisfies many properties of 
the usual convolution in $\RR^{d}$. The weighted spherical means, 
$T_\theta^\kappa$, with respect to $h_\kappa^2$ is defined implicitly by 
the formula 
\begin{equation} \label{eq:2.7}
 b_\lambda \int_0^\pi T_\theta^\kappa f(x) g(\cos \theta) 
(\sin \theta)^{2\lambda} d\theta = (f\star_\kappa g)(x), 
\qquad 0\le \theta \le \pi. 
\end{equation} 
where $g$ is any $L^1(w_\lambda)$ function. For $\kappa =0$, $V_\kappa = id$, 
the weighted spherical means coincide 
with the weighted means $T_\theta f$ in (\ref{eq:1.10}). Many properties 
of $T_\theta f$, given in \cite{BBP,P}, can be extended to the weighted 
means $T_\theta^\kappa f$. In particular, we have 
$$
 \|T_\theta^\kappa f\|_{\kappa,p} \le \|f\|_{\kappa,p} \qquad \hbox{and}
  \qquad \lim_{\theta \to 0} \|T_\theta^\kappa f -f \|_{\kappa,p} =0.
$$
Consequently, the following definition of a modulus of smoothness, 
$\omega_r(f;t)_{\kappa,p}$, makes sense. Let $r >0$, for $f \in 
L^p(h_\kappa^2)$, $1 \le p < \infty$, or $f \in C(S^{d-1})$, define
\begin{equation} 
\omega_r(f,t)_{\kappa,p} 
  := \sup_{0\le \theta \le t} \|(I-T_\theta^\kappa)^{r/2}\|_{\kappa,p}.
\label{eq:2.8}
\end{equation}
For the unweighted case ($\kappa =0$), such a definition was given in 
\cite{Rus} and the case $r$ being an even integer had appeared in several 
early references (see the discussion in \cite{Rus}). One of the important
properties of this modulus of smoothness is that it is equivalent to a
K-functional. 

Let $r > 0$. Recall the equation (\ref{eq:2.1}). We define 
$(-\Delta_{h,0})^{r/2}g$ by 
$$
(-\Delta_{h,0})^{r/2} g 
 \sim \sum_{n=1}^\infty (n(n+2\lambda_\kappa))^{r/2} Y_n(h_\kappa^2;g).
$$   
Furthermore, define the function space $\CW_r^p(h_\kappa^2)$ by
$$
\CW^p_r(h_\kappa^2) =\left \{f \in L^p(h_\kappa^2): 
 (k (k+2 \lambda))^{\frac{r}2} P_k(h_\kappa^2;f) = 
P_k(h_\kappa^2;g) \,\, \hbox{some} \,\, g \in   L^p(h_\kappa^2)\right \}.
$$
The K-functional between $L^p(h_\kappa^2)$ and $\CW_r^p(h_\kappa^2)$ is 
defined by 
\begin{equation}\label{eq:K-funct}
K_r(f;t)_{\kappa,p} := \inf \big \{ \|f - g\|_{\kappa,p} + t^r 
 \|(-\Delta_{h,0})^{r/2}\, g\|_{\kappa,p}, \; g \in \CW_r^p(h_\kappa^2) \big\}.
\end{equation}
It is equivalent to the modulus of smoothness in the following sense: 

\bigskip

\begin{theorem} \label{thm:2.1}
For $1 \le p \le \infty$, there exist two positive constants $c_1$ and $c_2$
such that for $f \in L^p(h_\kappa^2)$, 
$$
c_1 \omega_r(f;t)_{\kappa,p} \le K_r(f;t)_{\kappa,p} \le 
 c_2\,  \omega_r(f;t)_{\kappa,p}.
$$
\end{theorem}

The modulus of smoothness or the K-functional can be used to characterize
the best approximation by polynomials. Let 
$$
   E_n(f)_{\kappa,p} : = \inf \left \{ \|f - P \|_{\kappa,p}: 
        P \in \Pi_n^{d} \right\}, 
$$
We state the direct and the inverse theorems in terms of the modulus of 
smoothness.

\bigskip

\begin{theorem} \label{thm:2.2}
For $f \in L^p(h_\kappa^2)$, $1 \le p \le \infty$, 
$$
   E_n(f)_{\kappa,p} \le c\, \omega_r(f;n^{-1})_{\kappa,p}.
$$
On the other hand, 
$$
 \omega_r(f,n^{-1})_{\kappa,p} \le c\, n^{-r} \sum_{k=0}^n (k+1)^{r-1} 
    E_k(f)_{\kappa,p}.
$$
\end{theorem}

These theorems are proved in \cite{X02a,X02b}, following closely the method
developed in \cite{Rus} where these theorems were essentially established in 
the case $\kappa =0$. For results in the unweighted cases, see also
\cite{LW}. The problem of best approximation has been studied by many authors.
We refer to \cite{BBP,Kam,LN,P,Ra,Rus,W} and the references therein.

One can study summability of $h$-harmonic expansions (see \cite{X97b,X01a,LX})
and the convolution structure $\star_{\kappa}$ is useful in this direction
(see \cite{X02a,X03}). For summability of the ordinary harmonic expansions
(unweighted case), we refer to \cite{BBP,BL,BC,Kam,LW,LN,P,Rus} and the 
references therein.

In the case of usual surface measure on $S^{d-1}$, there are various other 
moduli of smoothness that have been used to characterize the best 
approximation; see, for example, \cite{Rus2,W}. It would be nice to define
one that can be given explicitly for the weighted case. 

\section{Generalized translation operator and Approximation on $B^d$}
\setcounter{equation}{0}

Recall the weight function $W_{\kappa,\mu}^B(x)$ defined in (\ref{eq:1.6}),
in which $h_\kappa$ is an reflection invariant weight function defined on 
$\RR^d$. Let $a_{\kappa,\mu}$ denote the normalization constant for 
$W_{\kappa,\mu}^B$. Denote by $L^p(W_{\kappa,\mu}^B)$, $1 \le p \le \infty$, 
the space of measurable functions defined on $B^d$ with the finite norm
$$
  \|f\|_{W_{\kappa,\mu}^B,p}:= \Big(a_{\kappa,\mu} \int_{B^d} |f(x)|^p 
     W_{\kappa,\mu}^B(x) dx \Big)^{1/p}, \qquad 1 \le p < \infty,
$$
and for $p = \infty$ we assume that $L^\infty$ is replaced by $C(B^d)$, the 
space of continuous function on $B^d$.

\subsection{Background}
Let $\CV_n^d(W_{\kappa,\mu}^B)$ denote the space of orthogonal polynomials of 
degree $n$ with respect to $W_{\kappa,\mu}^B$ on $B^d$. Elements of 
$\CV_n^d(W_{\kappa,\mu}^B)$ are closely related to the $h$-harmonics associated
with the weight function 
$$
h_{\kappa,\mu}(y_1,\ldots,y_{d+1}) = h_\kappa(y_1,\ldots,y_d)|y_{d+1}|^{\mu}
$$
on $\RR^{d+1}$, where $h_\kappa$ is associated with the reflection group $G$. 
The function $h_{\kappa,\mu}$ is invariant under the group $G \times \ZZ_2$. 
Let $Y_n$ be such an $h$-harmonic polynomial of degree $n$ and assume that 
$Y_n$ is even in the $(d+1)$-th variable; that is, $Y_n(x,x_{d+1})=Y_n(x,
-x_{d+1})$. We can write
\begin{equation} \label{eq:3.1}
   Y_n(y) = r^n P_n(x), \qquad  y = r(x,x_{d+1}) \in \RR^{d+1}, \quad 
 r = \|y\|, \quad (x,x_{d+1}) \in S^d,
\end{equation}
in polar coordinates. Then $P_n$ is an element of $\CV_n^d(W_{\kappa,\mu}^B)$  
and this relation is an one-to-one correspondence (\cite{X01a}). Under the 
changing variables $y \mapsto r(x,x_{d+1})$, $h_{\kappa,\mu}$ becomes 
$W_{\kappa,\mu}^B$ and the elementary formula 
\begin{equation} \label{eq:3.2}
   \int_{S^d} P(y) d\omega = \int_{B^d} \left[ P(x,\sqrt{1-\|x\|^2}\,)+
     P(x,-\sqrt{1-\|x\|^2}\,) \right]\frac{dx}{\sqrt{1-\|x\|^2}}
\end{equation}
shows the relation between their normalization constants. 

Let $\Delta_h^{\kappa,\mu}$ denote the $h$-Laplacian associated with 
$h_{\kappa,\mu}$ and $\Delta_{h,0}^{\kappa,\mu}$ denote the corresponding
spherical $h$-Laplacian. When $\Delta_h^{\kappa,\mu}$ is applied to functions 
on $\RR^{d+1}$ that are even in the $(d+1)$-th variable, the spherical 
$h$-Laplacian can be written in polar coordinates $y = r(x,x_{d+1})$ as 
 (\cite{X01b}): 
$$
 \Delta_{h,0}^{\kappa,\mu} = \Delta_h - \langle x, \nabla 
 \rangle^2 - 2 (\lambda_\kappa+\mu) \langle x, \nabla \rangle, 
$$
in which the operators $\Delta_h$ and $\nabla=(\partial_1,\ldots,\partial_d)$ 
are all acting on $x$ variables and $\Delta_h$ is the $h$-Laplacian associated
with $h_\kappa$ on $\RR^d$. Define 
\begin{equation*} 
 D_{\kappa,\mu}^B := \Delta_h - \langle x, \nabla \rangle^2 - 
        2 (\lambda_\kappa+\mu) \langle x, \nabla \rangle, 
\end{equation*}
as an operator acting on functions defined on $B^d$. It follows that the 
elements of $\CV_n^d(W_{\kappa,\mu})$ are eigenfunctions of $D_{\kappa,\mu}^B$:
\begin{equation} \label{eq:3.3}
  D_{\kappa,\mu}^B P = - n (n+2\lambda_\kappa+2\mu) P,    
    \qquad P \in\CV_n^d(W_{\kappa,\mu}^B).
\end{equation}
For the classical weight function $W_\mu^B(x) = (1-\|x\|^2)^{\mu-1/2}$, the 
operator $D_{\kappa,\mu}^B$ becomes a pure differential operator which 
is a classical result going back to Hermite (see \cite{AF} and 
\cite[Chapt. 12]{Er}). 

For $f \in L^2(W_{\kappa,\mu}^B)$, its orthogonal expansion is given by
$$
 L^2(W_{\kappa,\mu}^B) = \sum_{n=0}^\infty\bigoplus \CV_n^d(W_{\kappa,\mu}^B)
 : \qquad f = \sum_{n=0}^\infty \proj_n^{\kappa,\mu} f, 
$$
where $\proj_n^{\kappa,\mu}:L^2(W_{\kappa,\mu}^B)\mapsto
\CV_n^d(W_{\kappa,\mu}^B)$ 
is the projection operator, which can be written as an integral
\begin{equation} \label{eq:3.4}
  \proj_n^{\kappa,\mu} f(x) = a_{\kappa,\mu} \int_{B^d} f(y) 
      P_n(W_{\kappa,\mu}^B; x,y) W_{\kappa,\mu}^B(y)dy, 
\end{equation}
where $P_n(W_{\kappa,\mu}^B; x,y)$ is the reproducing kernel of 
$\CV_n^d(W_{\kappa,\mu}^B)$.  
The intertwining operator associated with $h_{\kappa,\mu}$, denoted by 
$V_{\kappa,\mu}$, is given in terms of the intertwining operator $V_\kappa$ 
associated to $h_\kappa$ and the operator $V_\mu^{\ZZ_2}$ associated to 
$h_\mu(x) =|x_{d+1}|^\mu$, $x \in \RR^{d+1}$, which is given explicitly by
(\ref{eq:2.5}) (setting $d=1$ and $\kappa_1=\mu$ there); that is,
$$
  V_{\kappa,\mu} f(x,x_{d+1}) = c_\mu \int_{-1}^1 
         V_\kappa [f(\cdot,x_{d+1}t)](x) (1+t)(1-t^2)^{\mu-1}dt,
$$
where $x \in \RR^d$. Since polynomials in $\CV_n(W_{\kappa,\mu}^B)$ correspond
to $h$-harmonics that are even in the last coordinates, we introduce a modified
operator
\begin{align} \label{eq:3.5}
  V_{\kappa,\mu}^B f(x,x_{d+1}) & :=  \left[V_{\kappa,\mu} f(x,x_{d+1})
   + V_{\kappa,\mu} f(x, -x_{d+1}) \right] \\
  & =  c_\mu \int_{-1}^1 
         V_\kappa [f(\cdot,x_{d+1}t)](x) (1-t^2)^{\mu-1}dt, \notag
\end{align}
acting on functions defined on $\RR^{d+1}$. 

\subsection{Generalized translation operator} For weight $W_{\kappa,\mu}^B$ 
on $B^d$, we define a convolution, denoted by $\star_{\kappa,\mu}^B$, as 
follows: For $f\in L^1(W_{\kappa,\mu}^B)$ and 
$g \in L^1(w_{\lambda_\kappa+\mu},[-1,1])$, 
$$
  (f \star_{\kappa,\mu}^B g)(x) = a_{\kappa,\mu} \int_{B^d} f(y) 
     V_{\kappa,\mu}^Bg(\langle X, \cdot \rangle)(Y)W_{\kappa,\mu}^B(y)dy,
$$
where $X = (x, \sqrt{1-\|x\|^2})$ and $Y = (y, \sqrt{1-\|y\|^2})$. The 
properties of this convolution can be derived from the corresponding 
convolution on the sphere. Let $f \star_{\kappa,\mu} g$ denote the convolution
defined in (\ref{eq:2.6}) with respect to $h_{\kappa,\mu}$. In fact, 
(\ref{eq:3.2}) immediately implies the following proposition.

\bigskip

\begin{proposition}  \label{prop:3.1}
For $f\in L^1(W_{\kappa,\mu}^B)$ and $g \in L^1(w_{\lambda_\kappa+\mu},
[-1,1])$, 
$$
(f \star_{\kappa,\mu}^B g)(x) = (F\star_{\kappa,\mu} g)(x, \sqrt{1-\|x\|^2}),
 \qquad \hbox{where} \qquad F(x,x_{d+1}) := f(x).
$$
\end{proposition}  

We now define the generalized translation operator on $B^d$ implicitly via 
the convolution operator. 

\bigskip

\begin{definition}\label{defn3}
For $f\in L^1(W_{\kappa,\mu}^B)$, the generalized translation operator 
\allowbreak
$T_\theta(W_{\kappa,\mu}^B;f)$ is defined implicitly by 
\begin{align} \label{eq:3.6}
b_{\lambda+\mu} \int_0^\pi T_\theta(W_{\kappa,\mu}^B;f,x)g(\cos \theta) 
   (\sin \theta)^{2\lambda+2\mu}d\theta 
 = (f\star_{\kappa,\mu}^B g)(x),
\end{align}
where $\lambda = \lambda_\kappa$, for every $g \in L^1(w_{\lambda+\mu},
[-1,1])$. 
\end{definition}

The generalized translation operator $T_\theta(W_{\kappa,\mu}^B)$ is related
to the weighted spherical means associated with the weight function 
$h_{\kappa,\mu}$ on $S^d$. For $F \in L^1(h_{\kappa,\mu}^2)$, denote the
weighted spherical means by $T_\theta^{\kappa,\mu} F$ as defined in 
(\ref{eq:2.7}). From the definitions of $T_\theta(W_{\kappa,\mu}^B)$ and 
$T_\theta^{\kappa,\mu}$, Proposition \ref{prop:3.1} shows that the following
relation holds:

\bigskip

\begin{proposition} \label{prop:3.3}
For each $x \in B^d$ the operator $T_\theta(W_{\kappa,\mu}^B;f,x)$ is a 
uniquely determined $L^\infty$ function in $\theta$. Furthermore, define 
$F(x,x_{d+1}) = f(x)$; then
\begin{equation} \label{eq:3.7}
 T_\theta(W_{\kappa,\mu}^B;f,x) = T_\theta^{\kappa,\mu} F(x,x_{d+1}), \qquad
  x \in B^d, \quad x_{d+1} = \sqrt{1-\|x\|^2}.
\end{equation} 
\end{proposition}

We could of course define the generalized translation operator by the formula
(\ref{eq:3.7}). The convolution $\star_{\kappa,\mu}^B$, however, will be used
in the following section. These relations allow us to derive the following 
properties of the generalized translation operator.

\bigskip

\begin{proposition} \label{prop:3.4}
The generalized translation operator $T_\theta(W_{\kappa,\mu}^B; f)$ satisfies
the following properties:
\begin{enumerate}
\item Let $f_0(x) =1$, then $T_\theta(W_{\kappa,\mu}^B; f_0,x) =1$; 
\item For $f\in L^1(W_{\kappa,\mu}^B)$, 
$$
  \proj_n^{\kappa,\mu} T_\theta(W_{\kappa,\mu}^B; f) 
      = \frac{C_n^{\lambda_\kappa+\mu}(\cos \theta)}
    {C_n^{\lambda_\kappa+\mu}(1)} \proj_n^{\kappa,\mu} f; 
$$
\item $T_\theta(W_{\kappa,\mu}^B): \Pi_n^d \mapsto \Pi_n^d$ and
$$
  T_\theta(W_{\kappa,\mu}^B; f) \sim \sum_{n=0}^\infty 
   \frac{C_n^{\lambda_\kappa+\mu}(\cos \theta)}{C_n^{\lambda_\kappa+\mu}(1)}
      \proj_n^{\kappa,\mu}f;  
$$
\item For $0 \le \theta \le \pi$,
$$
 T_\theta(W_{\kappa,\mu}^B; f) - f = \int_0^\theta 
   (\sin s)^{-2\lambda_\kappa-2\mu} 
ds\int_0^s T_t(W_{\kappa,\mu}^B;D_{\kappa,\mu}^B f) 
   (\sin t)^{2\lambda_\kappa+2\mu}dt;
$$
\item For $f\in L^p(W_{\kappa,\mu}^B)$, $1 \le p < \infty$, or $f\in C(B^d)$,
$$
 \|T_\theta(W_{\kappa,\mu}^B; f)\|_{W_{\kappa,\mu}^B,p} \le 
   \|f\|_{W_{\kappa,\mu}^B,p} \quad \hbox{and}  \quad 
\lim_{\theta \to 0}\|T_\theta(W_{\kappa,\mu}^B; f)-f\|_{W_{\kappa,\mu}^B,p}=0.
$$
\end{enumerate}
\end{proposition}

\begin{proof}
All these properties follow from the integral relation (\ref{eq:3.2}), the 
relation 
\begin{align}\label{eq:3.8}
P_n(W_{\kappa,\mu}^B; x,y) = \frac{1}{2} \Big[
 & Y_n(h_{\kappa,\mu}^2; (x,\sqrt{1-\|x\|^2}),(y,\sqrt{1-\|y\|^2}))\\
 & + Y_n(h_{\kappa,\mu}^2; (x,\sqrt{1-\|x\|^2}),(y,-\sqrt{1-\|y\|^2})) \Big], 
  \notag
\end{align}
the connection between $T_\theta(W_{\kappa,\mu}^B; f)$ and 
$T_\theta^{\kappa,\mu}F$ in Proposition \ref{prop:3.4}, and the corresponding 
relations for $T_\theta^{\kappa,\mu} F$ in \cite{X02b}. Recall the projection 
operator $Y_n(h_{\kappa,\mu}^2;F)$ of $h$-harmonics defined in (\ref{eq:2.2}).
The relations (\ref{eq:3.2}) and (\ref{eq:3.8}) show that 
$$
  \proj_n^{\kappa,\mu} f(x) = Y_n(h_{\kappa,\mu}^2; F, X), \qquad 
    X=(x,\sqrt{1-\|x\|^2}).
$$
Hence, it follows from Proposition \ref{prop:3.4} and property (2) of 
Proposition 2.4 in \cite{X02b} that 
\begin{align*}
 \proj_n^{\kappa,\mu} T_\theta(W_{\kappa,\mu}^B;f,x) & = 
   \proj_n^{\kappa,\mu} T_\theta^{\kappa,\mu} F(X) = 
   Y_n(h_{\kappa,\mu}^2; T_\theta^{\kappa,\mu}F, X)\\
& = \frac{C_n^{\lambda + \mu}(\cos\theta)}{C_n^{\lambda + \mu}(1)} 
  Y_n(h_{\kappa,\mu}^2; F, X)
 =  \frac{C_n^{\lambda + \mu}(\cos\theta)}{C_n^{\lambda + \mu}(1)} 
    \proj_n^{\kappa,\mu} f(x).
\end{align*}
This proves (2). The property (3) follows from (2). Since $F(x,x_{d+1})=f(x)$ 
is evidently even in $x_{d+1}$, the definition shows that 
$\Delta_0^{\kappa,\mu} F(x,x_{d+1}) = D_{\kappa,\mu}^B f(x)$. Consequently,
by the definition of $T_\theta^{\kappa,\mu}$ and (\ref{eq:3.7}), we conclude
that 
$$
T_\theta^{\kappa,\mu} \Delta_0^{\kappa,\mu} F(x,x_{d+1}) = 
 T_\theta(W_{\kappa,\mu}^B; D_{\kappa,\mu}^B f, x),
$$
from which the property (4) follows from the corresponding property of 
$T_\theta^{\kappa,\mu}$ (Proposition 2.4 in \cite{X02b}). Finally, to prove 
the property (5), we note that it follows from the definition of 
$V_{\kappa,\mu}$ at (\ref{eq:3.5}) that 
$$
V_{\kappa,\mu}\left[ g(\langle (x,-x_{d+1}),\cdot\rangle)\right](y,y_{d+1}) 
= V_{\kappa,\mu}\left[ g(\langle (x, x_{d+1}),\cdot\rangle)\right](y,-y_{d+1})
$$
Hence, the definition of $T_\theta^{\kappa,\mu}$ shows that 
$T_\theta^{\kappa,\mu} F (x, x_{d+1})= T_\theta^{\kappa,\mu} F (x, -x_{d+1})$.
Consequently, by (\ref{eq:3.2}),
\begin{align*}
\int_{B^d} \left|T_\theta(W_{\kappa,\mu}^B; f, x)\right|^p W_{\kappa,\mu}^p(x)
 dx & = 
\int_{B^d} \left|T_\theta^{\kappa,\mu} F(x,\sqrt{1-\|x\|^2})\right|^p 
  W_{\kappa,\mu}^p(x)dx \\
 & =  \int_{S^d} \left|T_\theta^{\kappa,\mu} F(y)\right|^p
    h_{\kappa,\mu}^2(y) dy. 
\end{align*}
Let $\|\cdot\|_{\kappa,\mu,p}$ denote the $L^p(h_{\kappa,\mu}^2)$ norm on
$S^d$. For $f(x) = F(X)$, we have $\|f\|_{W_{\kappa,\mu}^B,p} = 
\|F\|_{\kappa,\mu,p}$. Hence, the property (5) follows from the corresponding 
property of $T_\theta^{\kappa,\mu}$ (Proposition 2.4 in \cite{X02b}).
\end{proof} 

Recall the definition of the generalized translation operator for the
Gegenbauer expansions at (\ref{eq:1.1}). The reason that 
$T_\theta(W_{\kappa,\mu}^B)$ is called the generalized translation operator 
lies in the property (2), since for $d=1$ and $\kappa =0$ the property (2) 
agrees with (\ref{eq:1.3}). 
 
Once the generalized translation operator is defined, we see that 
(\ref{eq:3.6}) expresses the convolution of $f\star_{\kappa,\mu}^B g$ as an
integral of one variable. For $g\in L^1(w_\lambda,[-1,1])$, its Gegenbauer
expansion can be written as 
$$
 g(t) \sim \sum_{n=0}^\infty \wh g_n^\lambda \frac{n+\lambda}{\lambda}
   C_n^\lambda(t), \qquad \hbox{where} \qquad 
 \wh g_n^\lambda = b_\lambda \int_{-1}^1 g(s) 
     \frac{C_n^\lambda (s)}{C_n^\lambda (1)} w_\lambda(s)ds,
$$
since the $L^2(w_\lambda,[-1,1])$ norm of $C_n^\lambda$ is equal to 
$C_n^\lambda(1) \lambda /(n+\lambda)$. Hence, it follows from the property
(2) of Proposition \ref{prop:3.4} that 
$$
 \proj_n^{\kappa,\mu} (f \star_{\kappa,\mu}^B g) = \wh g_n^{\lambda+\mu} 
     \proj_n^{\kappa,\mu} f,  
$$
which is the analog of the familiar property $\wh {f\ast g} = \wh f \cdot 
\wh g$ of the ordinary convolution. The convolution $\star_{\kappa,\mu}^B$ also
satisfies several other properties of the ordinary convolution. For example, 
it satisfies Young's inequality: 

\bigskip

\begin{proposition}
For $f \in L^q(W_{\kappa,\mu}^B)$ and $g \in L^r(w_{\lambda+\mu};[-1,1])$, 
$$
\|f\star_{\kappa,\mu}^B g\|_{W_{\kappa,\mu}^B,p}\le \|f\|_{W_{\kappa,\mu}^B,q} 
   \|g\|_{w_{\lambda+\mu},r}, 
$$
where $p,q,r \ge 1$, $p^{-1} = r^{-1}+q^{-1}-1$ and 
$\|\cdot\|_{w_{\lambda+\mu},r}$ denotes the $L^r(w_{\lambda+\mu},[-1,1])$ norm.
\end{proposition}

This follows from Proposition \ref{prop:3.1} and Young's inequality for
$\star_{\kappa,\mu}$ in \cite{X02b}.

\subsection
{Generalized translation operator for the classical weight $W_\mu^B$}
Recall the integral formula \eqref{eq:1.1} of the generalized translation
operator for the weight function $w_\lambda$. In the case of the classical 
weight function $W_\mu^B$ in \eqref{eq:1.8}, it is possible to give an 
integral formula for the generalized translation operator in the same spirit.

To see how such a formula may look like, we turn to a further relation between
functions on $B^d$ and those on $S^{d+m}$, where $m$ is a positive integer.
For $f(x)$ on $B^d$, define $F(x,x') = f(x)$ on $\RR^{d+m}$. Then
$$
  \int_{S^{d+m}} F(y) d\omega(y) = \int_{B^d} f(x) (1-\|x\|^2)^{(m-1)/2} dx. 
$$ 
As shown in \cite{X01a}, this relation preserves orthogonal structure. This 
suggests a relation between $T_\theta(W_\mu^B;f)$ on $B^d$ and the ordinary
spherical means $T_\theta f$ on $S^{d+m}$, similar to the one in Proposition 
\ref{prop:3.3} (which is the case $m=0$). On the other hand, the spherical
means is given by the formula \eqref{eq:1.10}. Hence, it is possible to 
derive a formula for $T_\theta (W_\mu^B;f)$ from that of spherical means.
It is this heuristic argument that suggests the following formula. 

Let $I$ denote the $d\times d$ identity matrix and define the symmetric matrix
$$
  A(x) = (1-\|x\|^2) I - x^T x, \qquad x = (x_1,\ldots, x_d),
$$
where $x^T$ is the transpose of $x$ ($x^T$ is a column vector). In the 
following, we take $u\in \RR^d$ as a row vector. For $u \in \RR^d$, the 
inequality $1- uA(x)u^T \ge 0$ defines an ellipsoid in $\RR^d$ (see below). 

\bigskip

\begin{theorem} \label{thm:3.5}
For $W_\mu^B$ in \eqref{eq:1.8}, the generalized translation operator
is given by 
\begin{align} \label{eq:3.9}
& T_\theta(W_\mu^B;f,x) = A_\mu \left(\sqrt{1-\|x\|^2}\right)^{d-1} \\ 
 & \times  \int_{uA(x)u^T \le 1} f\left(\cos \theta x + \sin \theta \sqrt{1-\|x\|^2}\,
     u\right)(1 - u A(x) u^T)^{\mu-1} du,\notag
\end{align}
where $A_\mu = 1/ \int_{B^d} (1-\|x\|^2)^{\mu-1} dx$ is the normalization 
constant for $W_{\mu -1/2}^B$. 
\end{theorem}

\begin{proof}
Although the explicit formula of $V_\mu^B$ is known for $W_\mu^B$, it does
not seem to be easy to verify the defining formula \eqref{eq:3.6} directly. 
Instead, we will verify the property (2) of Proposition \ref{prop:3.4}. In 
other words, let $T_\theta^* f$ denote the right hand side of \eqref{eq:3.9}; 
we show that $T_\theta^*f = f$ for all $f \in \CV_n^d(W_\mu^B)$.
It is known that one basis of $\CV_n^d(W_\mu^B)$ consists of functions of the 
form $C_n^{\mu+(d-1)/2}(\langle x,y\rangle)$, $y \in S^d$ (\cite{X00}). Hence, 
it is sufficient to show that 
$$
  T_\theta^* C_n^{\mu+(d-1)/2}(\langle x,y\rangle) =
   \frac{C_n^{\mu+(d-1)/2}(\cos \theta)}{C_n^{\mu+(d-1)/2}(1)}
   C_n^{\mu+(d-1)/2}(\langle x,y\rangle), \qquad y \in S^d. 
$$
The matrix $A(x)$ has eigenvalues $1$ and $\sqrt{1-\|x\|^2}$ (repeated $d-1$ 
times) and it is symmetric. Hence, there is a unitary matrix $U$ such that 
$$
A(x) = U(x) \Lambda(x) U(x)^T, \qquad 
\Lambda(x) = \diag\left\{1,\sqrt{1-\|x\|^2}, \ldots,\sqrt{1-\|x\|^2}\right\}.
$$
The columns of $U(x)$ are the eigenvalues of $A(x)$. In particular, the
first column of $A(x)$ is $x/\|x\|$ and the other columns of $U(x)$ form an
orthonormal basis of the null space of $x^T x$; that is, the other columns
are mutually orthonormal and are also orthogonal to $x$. Changing variables 
$u \mapsto u U(x) := v$, the quadratic form becomes 
$$
  u A(x)u^T = v D(x) v^T = v_1^2 + \sqrt{1-\|x\|^2}( v_2^2 + \ldots + v_d^2), 
$$
which suggests one more change of variables $v \mapsto \sqrt{1-\|x\|^2} \,v
 D^{-1}(x) : = s$ with 
$$
  D(x) = \diag\left\{ \sqrt{1-\|x\|^2}, 1,\ldots 1\right\},
$$ 
so that the quadratic form becomes $u A(x)u^T = s s^T$. Hence, the integral
domain $u A(x) u^T \le 1$ becomes $B^d$ in $s$ variables. Since $U(x)$ is
unitary, we have $du = dv = ds / (1-\|x\|^2)^{(d-1)/2}$. Consequently, we have
\begin{align*}
 & T_\theta^* C_n^{\mu+(d-1)/2}(\langle x,y\rangle) \\
 & = a_\kappa 
    \int_{B^d} C_n^{\mu+(d-1)/2}\left(\cos \theta \langle x,y\rangle
    + \sin \theta \langle s,y U(x) D(x)\rangle \right) (1-\|s\|^2)^{\mu-1} ds,
\end{align*}
where we have used the fact that $\langle s D(x)U^T(x),y\rangle
= \langle s,y U(x) D(x)\rangle$. Since the first column of $U(x)$ is $x /\|x\|$
and $U$ is unitary, the vector $y U(x) D(x)$ has norm
\begin{align*}
  \|yU(x)D(x)\|^2 = yU(x)D^2(x)U^T(x) y & = 
    yU(x)D^2(x)U^T(x) y \cr
  & = y y^T - yU(I - D^2)U^T y^T =  1-\langle x,y\rangle^2,   
\end{align*}
as $\|y\| =1$ and $I - D^2 = \diag\{\|x\|^2, 0,\ldots,0\}$. Hence, using the 
formula
$$
 A_\mu \int_{B^d} f(\langle x, y \rangle) (1-\|x\|^2)^{\mu-1}dx =
   b_{\mu+(d-3)/2} \int_{-1}^1 f(t \|y\|)(1-t^2)^{\mu+(d-3)/2} dt, 
$$
which can be easily verified as the left hand side is invariant under the
rotation, we conclude that 
\begin{align*}
 & T_\theta^* C_n^{\mu+(d-1)/2}(\langle x,y\rangle)  =  b_{\mu+(d-3)/2}\\
  & \times \int_{-1}^1 
  C_n^{\mu+(d-1)/2}\left(\cos \theta \langle x,y\rangle +
    \sin \theta \sqrt{1-\langle x,y\rangle^2} \, t \right)(1-t^2)^{\mu+(d-3)/2} dt.
\end{align*}
Using the product formula for the Gegenbauer polynomials finishes the proof.
\end{proof}

For $\mu = 0$, the integral formula of $T_\theta(W_\mu;f)$ holds under the 
limit $\mu \to 0$ and the integral domain becomes $u A(x) u^T =1$. This case
has been studied in \cite{G}. See also \cite{Rub} in which a generalized 
translation operator is defined for the weight function $x_{d+1}^\mu d\omega$
on $S_+^d = \{x \in S^d: x_{d+1} \ge 0\}$, which is related to 
$T_\theta(W_\mu;f)$, but no integral formula as above given there. 

If $d =1$, then $A(x) = 1$ and $1- u^TA(x)u = 1-u^2$. Hence, the formula
for $T_\theta(W_\mu^B;f,x)$ when $d =1$ becomes
$$
 T_\theta(W_\mu^B;f,x) =A_\mu \int_{|u| \le 1} 
f\left(\cos \theta x + \sin \theta \sqrt{1-x^2}\, u\right)(1 - u^2)^{\mu-1} du,
$$
which agrees with the formula of $T_{\cos\theta} f(x)$ in \eqref{eq:1.1}. 

The proof of the theorem also gives an alternative expression for 
$T_\theta(W_\mu^B;f)$.

\bigskip

\begin{corollary}
Let $U(x)$ be the unitary matrix whose first column is $x / \|x\|$ and 
$D(x)=\diag\{\sqrt{1-\|x\|^2}, 1,\ldots,1\}$. For $W_\mu^B$ in \eqref{eq:1.8},
$$ 
T_\theta(W_\mu^B;f,x) = A_\mu 
  \int_{B^d} f\left(\cos \theta x + \sin \theta \,s D(x) U(x) \right) 
     (1 - \|s\|^2)^{\mu-1} ds,
$$
\end{corollary}

In the above formula we take $x$ and $s$ as row vectors in $\RR^d$. 
Recall that the first column of $U(x)$ is $x /\|x\|$ and the other 
columns of $U(x)$ are orthonormal vectors that are orthogonal to $x$
(that is, $\langle x, \xi\rangle =0$). Hence, we can write the formula
for $T_\theta(W_\mu^B;f)$ as an explicit integral over $B^d$. For example, 
in the case of $d=2$, 
$$
U(x)=\frac{1}{\|x\|} \left(\begin{matrix} x_1 & -x_2 \\ x_2 & x_1
 \end{matrix}\right) \qquad \hbox{and} \qquad
D(x) = \left(\begin{matrix} \sqrt{1-\|x\|^2} & 0 \\ 0 & 1\end{matrix}\right).
$$
In this case, it is more convenient to use the polar coordinates 
$x = r \cos \phi$ and $y=r\sin \phi$, where $r \ge 0$ and $0\le\phi\le 2\pi$. 

\bigskip

\begin{corollary}
For $d =2$, $x = r x'$ with $x' = (\cos \phi, \sin\phi)$,
\begin{align*}
&T_\theta(W_\mu^B;f,x) = \\ 
&\quad A_\mu \int_{B^2} f\left(r \cos \theta x' + \sqrt{1-r^2}\,
 \sin \theta s_1 x' + \sin\theta s_2 x' \right) (1 - \|s\|^2)^{\mu-1} ds.
\end{align*}
\end{corollary}

For $d > 2$ the formula of $U(x)$ can be messy. For example, in the case $d$ 
is odd, we do not have a simple formula. On the other hand, there are simple
expressions of $U(x)$ for $d = 4, 8, \ldots$. As examples, we give the 
formula of $U(x)$ for $d =3$ and $d= 4$ below.
$$
\left(\begin{matrix} \frac{x_1}{\|x\|} & \frac{x_2}{\sqrt{x_1^2+x_2^2}}&
  \frac{x_3}{\|x\|} \\ 
\frac{x_2}{\|x\|} & \frac{-x_1}{\sqrt{x_1^2+x_2^2}}&
  -\frac{x_2 x_3}{x_1 \|x\|} \\ 
\frac{x_3}{\|x\|} & 0&
  \frac{x_2^2+x_3^2}{x_1\|x\|}  
\end{matrix}\right)
\qquad \hbox{and} \qquad
\frac{1}{\|x\|} \left(\begin{matrix} x_1 & x_2 &x_3&x_4 \\
  x_2 & -x_1 &-x_4 &x_3 \\x_3 & x_4 &-x_1& -x_2 \\x_4 & -x_3 &x_2&-x_1 
 \end{matrix}\right).
$$ 
Using these one can write down the formula of $T_\theta(W_\mu;f)$
as an explicit integral over $B^d$. 

An interesting problem is to find an integral expression for 
$T_\theta(W_{\kappa,\mu}^B;f)$ with respect to other weight functions. 
One should perhaps start with the case that $h_\kappa$ is given by the
product weight function \eqref{eq:1.5}. 

\subsection
{Modulus of smoothness, K-functional and best approximation}
The property (5) of Proposition \ref{prop:3.4} shows that the following 
definition of the modulus of smoothness makes sense:

\bigskip

\begin{definition}
Let $r >0$. Define 
$$
(I - T_\theta^\kappa)^{r/2} \sim \sum_{n=0}^\infty 
 \left(1-{C_n^{\lambda_\kappa+\mu}(\cos\theta)}/{C_n^{\lambda_\kappa+\mu}(1)}
      \right)^{r/2} \proj_n^{\kappa,\mu}f.
$$
For $f\in L^p(W_{\kappa,\mu}^B)$, $1\le p <\infty$, or $f\in C(B^d)$, define
$$
\omega(f;t)_{W_{\kappa,\mu}^B,p} : = \sup_{\theta \le t} 
  \left \|\left(I-T_\theta(W_{\kappa,\mu}^B)\right)^{r/2} 
f \right\|_{W_{\kappa,\mu}^B,p}.
$$
\end{definition}

Because of Proposition \ref{prop:3.4}, this modulus of smoothness is related 
to the modulus $\omega_r(f;t)_{\kappa,\mu,p}$, defined in (\ref{eq:2.8}) but
associated with $h_{\kappa,\mu}$. In fact, we have
\begin{equation}\label{eq:3.10}
\omega(f;t)_{W_{\kappa,\mu}^B,p} = \omega_r(F;t)_{\kappa,\mu,p},
\qquad F(x,x_{d+1}) = f(x).    
\end{equation}
Consequently, properties of $\omega_r(f;t)_{W_{\kappa,\mu}^B,p}$ can be 
easily obtained from those of $\omega_r(f;t)_{\kappa,\mu,p}$ (see 
Proposition 3.6 of \cite{X02b}).

\bigskip

\begin{proposition} \label{prop:3.6}
The modulus of smoothness $\omega_r(f;t)_{W_{\kappa,\mu}^B,p}$ satisfies:
\begin{enumerate}
\item $\omega_r(f;t)_{W_{\kappa,\mu}^B,p} \to 0$ if $t \to 0$;
\item $\omega_r(f;t)_{W_{\kappa,\mu}^B,p}$ is monotone nondecreasing on 
$(0,\pi)$;
\item $\omega_r(f+g,t)_{W_{\kappa,\mu}^B,p} \le 
  \omega_r(f;t)_{W_{\kappa,\mu}^B,p}+ \omega_r(g,t)_{W_{\kappa,\mu}^B,p}$; 
\item For $0 < s < r$, 
$$
\omega_r(f;t)_{W_{\kappa,\mu}^B,p} \le 2^{[(r-s+1)/2]} 
       \omega_s(f;t)_{W_{\kappa,\mu}^B,p};
$$
\item If $(-D_{\kappa,\mu}^B)^k f \in L^p(W_{\kappa,\mu}^B)$, $k \in \NN$, 
then for $r > 2 k$
$$
\omega_r(f;t)_{W_{\kappa,\mu}^B,p} 
 \le c\, t^{2k} \omega_{r-2k}((-D_{\kappa,\mu}^B)^k f;t)_{W_{\kappa,mu}^B,p}.
$$
\end{enumerate}
\end{proposition}

To justify the definition of this modulus of smoothness, we show that it is 
equivalent to a K-functional, defined using the differential-difference 
operator associated with $W_{\kappa,\mu}^B$ (see \eqref{eq:3.3}). Let 
\begin{align*}
  \CW_r^p(W_{\kappa,\mu}^B):=\left \{f \in L^p(W_{\kappa,\mu}^B): 
     (-D_{\kappa,\mu}^B)^{r/2}f \in L^p(W_{\kappa,\mu}^B)\right \},
\end{align*}
where the fractional power of $D_{\kappa,\mu}^B$ on $f$ is defined by 
$$
(-D_{\kappa,\mu}^B)^{r/2} f \sim  \sum_{n=0}^\infty 
   (n (n+2\lambda_\kappa+2\mu))^{r/2}
   \proj_n^{\kappa,\mu} f, \quad  f \in L^p(W_{\kappa,\mu}^B).
$$
The K-functional between $L^p(W_{\kappa,\mu}^B)$ and 
$\CW_r^p(W_{\kappa,\mu}^B)$ is defined by 
$$
 K_r(f;t)_{W_{\kappa,\mu}^B,p}:= \inf \left \{ \|f-g\|_{W_{\kappa,\mu}^B,p} + 
  t^r \|(-D_{\kappa,\mu}^B)^{r/2} g\|_{W_{\kappa,\mu}^B,p} \right \},
$$
where the infimum is taken over all $g \in \CW_r^p(W_{\kappa,\mu}^B)$.

\bigskip

\begin{theorem} \label{thm:3.7}
For $f \in L^p(W_{\kappa,\mu}^B)$, $1 \le p \le \infty$, 
$$
  c_1 \omega_r(f;t)_{W_{\kappa,\mu}^B,p} \le K_r(f;t)_{W_{\kappa,\mu}^B,p} 
    \le c_2 \omega_r(f;t)_{W_{\kappa,\mu}^B,p}.  
$$
\end{theorem}

\begin{proof}
Again let $F(x,x_{d+1}) = f(x)$. Denote by $K_r(F; t)_{\kappa,\mu,p}$ the
K-functional defined in (\ref{eq:K-funct}) but with respect to the weight
function $h_{\kappa,\mu}^2$. Because of \eqref{eq:3.10} and the 
equivalence between $K_r(F; t)_{\kappa,\mu,p}$ and 
$\omega_r(F;t)_{\kappa,\mu,r}$ (see Theorem \ref{thm:2.1}), we only need 
to show that $K_r(f;t)_{W_{\kappa,\mu}^B,p} = K_r(F; t)_{\kappa,\mu,p}$.

It follows from \eqref{eq:3.2} that $\|\Delta_{h,0}^{\kappa,\mu}F\|_{\kappa,
\mu,p} = \|D_{\kappa,\mu}^Bf\|_{W_{\kappa,\mu}^B,p}$. Hence, 
$$
  K_r(f;t)_{W_{\kappa,\mu}^B,p} = \inf \left\{\|F - g_e\|_{\kappa,\mu,p} + 
  t^r \|\Delta_{h,0}^{\kappa,\mu} g_e\|_{\kappa,\mu,p}\right\}
   := K_r^*(F; t)_{\kappa,\mu,p},
$$
where the infimum is taken over all $g_e(x,x_{d+1})\in
\CW_r^p(h_{\kappa,\mu}^2)$ that are even in $x_{d+1}$. Evidently, 
$K_r^*(F; t)_{\kappa,\mu,p} \ge  K_r(F; t)_{\kappa,\mu,p}$. To complete the 
proof we show that $K_r^* = K_r$. For any $\varepsilon > 0$, fix a $g \in 
\CW_r^p(W_{\kappa,\mu}^p)$ such that 
$$
K_r(F;t)_{\kappa,\mu,p} \ge \|F-g\|_{\kappa,\mu,p} + t^r 
    \|-\Delta_{h,0}^{\kappa,\mu} g\|_{\kappa,\mu,p} - \varepsilon.
$$
Since $h_{\kappa,\mu}$ corresponds to $G \times \ZZ_2$, the spherical Laplacian
$\Delta_{h,0}^{\kappa,\mu}$ commutes with the sign change in $x_{d+1}$. 
Consequently, setting $g_e(x,x_{d+1}) = [g(x,x_{d+1})+ g(x,-x_{d+1})]/2$,
so that $g_e$ is even in $x_{d+1}$, it follows that 
$ \|\Delta_{h,0}^{\kappa,\mu}g_e\|_{\kappa,\mu,p} \le 
  \|\Delta_{h,0}^{\kappa,\mu}g\|_{\kappa,\mu,p}$.  
This and the fact that $\|F-g_e\|_{\kappa,\mu,p} \le \|F-g\|_{\kappa,\mu,p}$,
as $F$ is even in $x_{d+1}$, show that $K_r^*(F; t)_{\kappa,\mu,p} \le 
K_r(F; t)_{\kappa,\mu,p}+\varepsilon$. As $\varepsilon > 0$ is arbitrary, 
the proof follows. 
\end{proof}

One immediate consequence of the above equivalence is the following property
of the modulus of smoothness, which does not follow trivially from the 
definition of $\omega_r(f; t)_{W_{\kappa,\mu}^B,p}$ but it is clear for
the K-functional.

\bigskip

\begin{corollary} \label{cor:3.8}
For $f \in L^p(W_{\kappa,\mu}^B)$, $1 \le p \le \infty$, 
$$
\omega_r(f,\delta t)_{W_{\kappa,\mu}^B,p}  
 \le  c\, \max\{1,\delta^r\} 
    \omega_r(f, t)_{W_{\kappa,\mu}^B,p}.
$$
\end{corollary}

The direct and the inverse theorems for the best approximation by polynomials 
in $L^p(W_{\kappa,\mu}^B)$ is characterized in \cite{X02b} by the K-functional.
The equivalence in Theorem \ref{thm:3.7} allows us to state the 
characterization in terms of the modulus of smoothness. 
 
\bigskip

\begin{theorem}
For $f \in L^p(W_{\kappa,\mu}^B)$, $1 \le p \le \infty$, 
$$
 E_n(f)_{W_{\kappa,\mu}^B,p} \le c\, \omega_r(f;n^{-1})_{W_{\kappa,\mu}^B,p}.
$$
On the other hand, 
$$
\omega_r(f;n^{-1})_{W_{\kappa,\mu}^B,p} \le c\, n^{-r}\sum_{k=0}^n (k+1)^{r-1} 
    E_k(f)_{W_{\kappa,\mu}^B,p}.
$$
\end{theorem} 

Finally, let us mention that, by \eqref{eq:3.4}, \eqref{eq:3.8} and 
\eqref{eq:2.3}, 
$$
 \proj_{\kappa,\mu} f = f \star_{\kappa,\mu}^B p_n, \qquad \hbox{where}\qquad
    p_n(t) = \frac{n+\lambda+\mu}{\lambda+\mu} C_n^{\lambda+\mu}(t).
$$
Hence, all summation methods of orthogonal expansions with respect to 
$W_{\kappa,\mu}^B$ can be written in the form of $f\star_{\kappa,\mu}^B g_r$, 
where $g_r$ is the same summation method applies to the Gegenbauer series 
evaluated at point $t=1$. Since $T_\theta(W_{\kappa,\mu}^B;f)$, thus 
$\omega_r(f;t)_{W_{\kappa,\mu}^B,p}$, is defined in terms of 
$f\star_{\kappa,\mu}^B g$, the modulus of smoothness is a convenient tool
for studying the summability of orthogonal expansions on $B^d$. For various
results on the summability of orthogonal expansions with respect to 
$W_{\kappa,\mu}^B$, see \cite{DX,LX,X01a,X02b,X03} and the references therein. 

\section{Generalized translation operator and Approximation on $T^d$}
\setcounter{equation}{0}

Recall the weight function $W_{\kappa,\mu}^T$ defined in \eqref{eq:1.7},
in which $h_\kappa$ is a weight function invariant under a reflection group 
$G_0$ and even in each of its variables. That is, $h_\kappa$
is invariant under the semi-product of a reflection group $G_0$ and the 
abilian group $\ZZ_2^d$. 

The definition of $L^p(W_{\kappa,\mu}^T)$, $1 \le p \le \infty$, is similar
to the case of $W_{\kappa,\mu}^B$. The notions such as the space of orthogonal
polynomials $\CV_n^d(W_{\kappa,\mu}^T)$ and the reproducing kernel 
$P_n(W_{\kappa,\mu}^T;x,y)$ are also defined similarly as in the case of 
$W_{\kappa,\mu}^B$. 
 
\subsection{Background}
Elements of $\CV_n^d(W_{\kappa,\mu}^T)$ are closely related to the orthogonal
polynomials in $\CV_{2n}^d(W_{\kappa,\mu}^B)$. Let us denote by $\psi$ the 
mapping 
$$ 
\psi: (x_1,\ldots,x_d) \in B^d \mapsto (x_1^2,\ldots,x_d^2) \in T^d
$$ 
and define $(f\circ \psi)(x_1,\ldots,x_d) = f(x_1^2,\ldots,x_d^2)$. The 
elementary integral 
\begin{equation}\label{eq:4.1}
  \int_{B^d} f(x_1^2,\ldots,x_d^2) dx = \int_{T^d} f(x_1,\ldots,x_d)
    \frac{dx}{\sqrt{x_1\cdots x_d}}.
\end{equation}
shows that $\|f\|_{W_{\kappa,\mu}^T,p}=\|f\circ \psi\|_{W_{\kappa,\mu}^B,p}$.
The mapping $R \mapsto P$ given by 
\begin{equation}\label{eq:4.2}
 P_{2n}(x) = (R_n \circ \psi)(x) \qquad x \in B^d
\end{equation} 
is a one-to-one mapping between $R_n \in \CV_n^d(W_{\kappa,\mu}^T)$ and  
$P_{2n} \in \CV_{2n}^d(W_{\kappa,\mu}^B; \ZZ_2^d)$, the subspace of polynomials
in $\CV_{2n}^d(W_{\kappa,\mu}^B)$ that are even in each of its variables 
(invariant under $\ZZ_2^d$). In particular, applying 
$D_{W_{\kappa,\mu}^B}$ on $P_{2n}$ leads to a second order 
differential-difference operator acting on $R_n$. We denote this operator 
by $D_{\kappa,\mu}^T$. Then (\cite{X01b})
\begin{equation} \label{eq:4.3}
  D_{\kappa,\mu}^T R = - n (n+\lambda_\kappa+\mu) R,    
    \qquad R \in\CV_n^d(W_{\kappa,\mu}^T), 
\end{equation}
For the  weight function \eqref{eq:1.9}, the operator is a second order 
differential operator, which takes the form 
$$
D_{\kappa,\mu}^T =
\sum_{i=1}^d x_i(1-x_i) \frac {\partial^2 P} {\partial x_i^2} - 
 2 \sum_{1 \le i < j \le d} x_i x_j \frac {\partial^2 P}{\partial x_i 
 \partial x_j} + \sum_{i=1}^d \left( \Big(\kappa_i +\frac{1}{2}\Big) -
    \lambda x_i \right) \frac {\partial P}{\partial x_i}    
$$
(recall $\mu = \kappa_{d+1}$ in this case). This is classical, already known
in \cite{AF} at least for $d=2$ (see also \cite[Chapt. 12]{Er}). 

In the following, we also denote by $\proj_n^{\kappa,\mu}: 
L^2(W_{\kappa,\mu}^T)
\mapsto \CV_n^d(W_{\kappa,\mu}^T)$ the orthogonal projection operator. For 
$f \in L^2(W_{\kappa,\mu}^T)$, it can be written as an integral
$$
  \proj_n^{\kappa,\mu} f(x) = a_{\kappa,\mu} \int_{T^d} f(y) 
      P_n(W_{\kappa,\mu}^T; x,y) W_{\kappa,\mu}^T(y)dy, 
$$
where $P_n(W_{\kappa,\mu}^T; x,y)$ is the reproducing kernel of 
$\CV_n^d(W_{\kappa,\mu}^T)$. The relation \eqref{eq:4.1} implies, in 
particular, that (\cite{X01b})
\begin{equation}\label{eq:4.4}
  P_n(W_{\kappa,\mu}^T;x,y) = \frac{1}{2^d} \sum_{\varepsilon \in \ZZ_2^d}
   P_{2n}\left(W_{\kappa,\mu}^B; x^{1/2}, \varepsilon y^{1/2}\right),
\end{equation}
where $x^{1/2} = (\sqrt{x_1}, \ldots, \sqrt{x_d})$ and $\varepsilon u = 
(\varepsilon_1 u_1,\ldots \varepsilon_d u_d)$. We define a useful operator, 
$V_{\kappa,\mu}^T$, acting on functions of $d+1$ variables, 
\begin{equation}\label{eq:4.5}
  V_{\kappa,\mu}^T F(x,x_{d+1}) = \frac{1}{2^d}\sum_{\varepsilon \in \ZZ_2^d}
           V_{\kappa,\mu}^B F (\varepsilon x, x_{d+1}).
\end{equation}     
The definition of $V_{\kappa,\mu}^T$ is justified by the following fact:
Let $p_n^{(\alpha,\beta)}(t)$ denote the orthonormal Jacobi polynomial of 
degree $n$ associated to the weight function $w_{\alpha,\beta} (t) = 
(1-t)^\alpha (1+t)^\beta$ on $t \in [-1,1]$. Using the relation
\begin{equation}\label{eq:4.6}
\frac{2 n+\lambda}{\lambda} C_{2n}^\lambda(t) = 
   p_n^{(\lambda-1/2,-1/2)}(1)p_n^{(\lambda-1/2,-1/2)}(2t^2-1), 
\end{equation}
we can write the reproducing kernel $P_n(W_{\kappa,\mu}^T;x,y)$ as
$$
   P_n(W_{\kappa,\mu}^T;x,y) = p_n^{(\lambda_\kappa+\mu-\frac12,-\frac12)}(1)
   V_{\kappa,\mu}^T  \left[p_n^{(\lambda_\kappa+\mu-\frac12,-\frac12)} 
         (2 \langle \cdot, Y^{1/2} \rangle^2-1)\right](X^{1/2}).
$$
where $X^{1/2} = \left(\sqrt{x_1},\ldots,\sqrt{x_d},\sqrt{1-|x|}\right)$. 

\subsection{Generalized translation operator}
The operator $ V_{\kappa,\mu}^T$ is used to define a convolution operator on 
$T^d$: 

\bigskip

\begin{definition}\label{defn:4.1}
For $f \in L^1(W_{\kappa,\mu}^T)$ and $g \in L^1(w_{\lambda_\kappa+\mu};
[-1,1])$, we define 
$$
(f \star_{\kappa,\mu}^T g)(x) = a_{\kappa,\mu} \int_{T^d} f(y)V_{\kappa,\mu}^T 
   \left[g \left(2 \langle X^{1/2},\cdot \rangle^2-1\right)\right]
   (Y^{1/2}) W_{\kappa,\mu}^T(y) dy. 
$$
\end{definition}

Recall that $|x| = x_1 + \ldots + x_d$. Using \eqref{eq:4.1}, it is not hard to
show that $f \star_{\kappa,\mu}^T g$ is related to the convolution structure 
$f \star_{\kappa,\mu}^B g$ on $B^d$ (\cite{X03}):

\bigskip

\begin{proposition}
For $f \in L^1(W_{\kappa,\mu}^T)$ and $g \in L^1(w_{\lambda+\mu};[-1,1])$, 
$$
  \left((f\star_{\kappa,\mu}^T g)\circ \psi \right)(x) = 
     \left((f\circ \psi)\star_{\kappa,\mu}^B g(2\{\cdot\}^2-1) \right)(x).
$$
\end{proposition} 

The generalized translation operator associated with $W_{\kappa,\mu}^T$ 
is again defined implicitly in terms of the convolution structure.

\bigskip

\begin{definition}\label{defn5}
For $f\in L^1(W_{\kappa,\mu}^T)$, the generalized translation operator 
\allowbreak $T_\theta(W_{\kappa,\mu}^T;f)$ is defined implicitly by 
\begin{align} \label{eq:4.7}
b_{\lambda+\mu} \int_0^\pi T_\theta(W_{\kappa,\mu}^T;f,x)g(\cos 2 \theta) 
   (\sin \theta)^{2\lambda+2\mu}d\theta = (f \star_{\kappa,\mu}^T g)(x),
\end{align}
where $\lambda = \lambda_\kappa$, for every $g \in L^1(w_{\lambda+\mu})$. 
\end{definition}

The definition of $V_{\kappa,\mu}^T$ ensures that the generalized translation 
operator on $T^d$ is related to the one on $B^d$. This also shows that the
operator $T_\theta(W_{\kappa,\mu}^T;f)$ is well-defined. 

\bigskip

\begin{proposition} \label{prop:4.4}
For each $x \in T^d$ the operator $T_\theta(W_{\kappa,\mu}^T;f,x)$ is a 
uniquely determined $L^\infty$ function in $\theta$. Furthermore, 
\begin{equation} \label{eq:4.8}
 \left(T_\theta(W_{\kappa,\mu}^T;f)\circ \psi\right) (x) = 
    T_\theta(W_{\kappa,\mu}^B;f\circ \psi, x), \qquad  x \in T^d. 
\end{equation} 
\end{proposition}

The proof follows from the definition and the elementary formula 
$\cos 2 \theta = 2 \cos^2 \theta -1$; see \cite{X03}. This relation allows
us to derive properties of $T_\theta(W_{\kappa,\mu}^T;f)$. 

\bigskip

\begin{proposition} \label{prop:4.5}
The generalized translation $T_\theta(W_{\kappa,\mu}^T; f)$ satisfies the 
following properties:
\begin{enumerate}
\item Let $f_0(x) =1$, then $T_\theta(W_{\kappa,\mu}^T; f_0,x) =1$. 
\item For $f\in L^1(W_{\kappa,\mu}^T)$, 
$$
  \proj_n^{\kappa,\mu} T_\theta(W_{\kappa,\mu}^T; f) 
      = \frac{P_n^{(\lambda_\kappa+\mu-1/2,-1/2)}(\cos 2\theta)}
    {P_n^{(\lambda_\kappa+\mu-1/2,-1/2)}(1)} 
         \proj_n^{\kappa,\mu} f. 
$$
\item $T_\theta(W_{\kappa,\mu}^T): \Pi_n^d \mapsto \Pi_n^d$ and
$$
  T_\theta(W_{\kappa,\mu}^T; f) \sim \sum_{n=0}^\infty 
    \frac{P_n^{(\lambda_\kappa+\mu-1/2,-1/2)}(\cos 2\theta)}
    {P_n^{(\lambda_\kappa+\mu-1/2,-1/2)}(1)}    
    \proj_n^{\kappa,\mu}f.  
$$
\item For $0 \le \theta \le \pi$,
$$
 T_\theta(W_{\kappa,\mu}^T; f) - f =2\int_0^\theta 
   (\sin s)^{-2\lambda_\kappa-2\mu} 
ds\int_0^s T_t(W_{\kappa,\mu}^T;D_{\kappa,\mu}^T f) 
   (\sin t)^{2\lambda_\kappa+2\mu}dt.
$$
\item For $f\in L^p(W_{\kappa,\mu}^T)$, $1 \le p < \infty$, or $f\in C(T^d)$,
$$
 \|T_\theta(W_{\kappa,\mu}^T; f)\|_{W_{\kappa,\mu}^T,p} \le 
   \|f\|_{W_{\kappa,\mu}^T,p} \quad \hbox{and}  \quad 
\lim_{\theta \to 0}\|T_\theta(W_{\kappa,\mu}^T; f)-f\|_{W_{\kappa,\mu}^T,p}=0.
$$
\end{enumerate}
\end{proposition}

\begin{proof}
The property (1) is an easy consequence of \eqref{eq:4.8} and the property
(1) of Proposition \ref{prop:3.4}. Let $f \in \CV_n(W_{\kappa,\mu}^T)$. Then
$f\circ \psi \in \CV_{2n}(W_{\kappa,\mu}^B)$. Hence, by Proposition 
\ref{prop:3.4}, 
$$
 T_\theta(W_{\kappa,\mu}^T;f,x_1^2,\ldots,x_d^2) = 
 T_\theta(W_{\kappa,\mu}^B;f\circ\psi,x) = 
\frac{C_{2n}^{\lambda_\kappa+\mu}(\cos\theta)}{C_{2n}^{\lambda_\kappa+\mu}(1)}
   (f \circ \psi)(x).
$$
This proves the properties (2) and (3) upon using the relation \eqref{eq:4.2}.
Let $f = \sum c_k R_k$, $R_k \in \CV_n(W_{\kappa,\mu}^T)$. By \eqref{eq:4.2},
$P = R\circ \psi \in \CV_n(W_{\kappa,\mu}^B)$. Then 
\begin{align} \label{eq:4.9}
(D_{\kappa,\mu}^T f) \circ\psi & = - \sum c_k k(k+\lambda_\kappa+\mu) R_k
 \circ \psi\\
& = - 2^{-1} \sum c_k  2k(2k+2\lambda_\kappa+2\mu) P_{2k} \notag \\
& =  2^{-1} \sum c_k D_{\kappa,\mu}^B P_{2k} 
= 2^{-1} D_{\kappa,\mu}^B (f\circ \psi), \notag
\end{align}
from which (4) follows from the property (4) of Proposition \ref{prop:3.4}. 
Finally, a change of variables $x \mapsto \psi(x)$ shows that
$$
 \|T_\theta(W_{\kappa,\mu}^T;f)\|_{W_{\kappa,\mu}^T,p} = 
     \|T_\theta(W_{\kappa,\mu}^T;f)\circ \psi\|_{W_{\kappa,\mu}^B,p}  
   =  \|T_\theta(W_{\kappa,\mu}^B;f\circ \psi)\|_{W_{\kappa,\mu}^B,p}, 
$$
which is less than or equal to $\|f\circ \psi\|_{W_{\kappa,\mu}^B,p} =
\|f\|_{W_{\kappa,\mu}^T,p}$ by the property (5) of Proposition \ref{prop:3.4}. 
\end{proof} 

Using the relation to $T_\theta(W_{\kappa,\mu}^B;f)$, we can derive from 
Theorem \ref{thm:3.5} an integral formula for $T_\theta(W_\mu^T;f)$, where 
$W_\mu^T(x) = (x_1\cdots x_d)^{-1/2}(1-|x|)^{\mu -1/2}$. One interesting
question is to find such a formula for the classical weight function
$W_\kappa^T$ in \eqref{eq:1.9}. 

In the case of $d =1$ and $G = \ZZ_2^d$, the weight function 
$W_{\kappa,\mu}^T$ becomes the Jacobi weight function $w_{\kappa,\mu}(t) 
= 2^{\kappa+\mu} t^\kappa (1-t)^\mu$ on $[0,1]$ (see \eqref{eq:1.9}), whose 
corresponding orthogonal polynomials are $P_n^{(\kappa,\mu)}(2t-1)$. 
The orthogonal expansion of $f$ in Jacobi polynomials is defined by 
$$
  f(t) \sim \sum_{n=0}^\infty a_n (f) p_n^{(\alpha,\beta)} (t), \quad
   \hbox{where} \quad
   a_n(f) = c_{\alpha,\beta} \int_{-1}^1 f(s) p_n^{(\alpha,\beta)}(s)ds 
$$
and  $c_{\alpha,\beta}^{-1} = \int_{-1}^1 w_{\alpha,\beta}(s)ds$. The 
usual generalized translation operator, $S_\theta f(t)$, for the Jacobi 
expansion is an operator defined by (\cite{AW})
\begin{equation*}
  S_\theta f(t) \sim \sum_{n=0}^\infty a_n (f) 
    p_n^{(\alpha,\beta)}(\cos \theta)
    p_n^{(\alpha,\beta)} (t). 
\end{equation*}
We should emphasis, however, that the operator $S_\theta f$ is different 
from the case $d =1$ of the generalized translation operator $T_\theta
(W_{\kappa,\mu}^T;f)$. Even in the case of $\alpha = \lambda + \mu -1/2$
and $\beta = -1/2$, they are different as can be seen from Proposition 
\ref{eq:4.5}.
The convolution structure of the Jacobi expansions defined via $S_\theta$ 
has a natural extension to the product Jacobi weight functions on the unit 
cube $[-1,1]^d$. The convolution structure defined above works for 
$W_{\kappa,\mu}^T$ on the simplex.   

\subsection
{Modulus of smoothness, K-functional and best approximation}
We can also define a modulus of smoothness on $T^d$ using the generalized
translation operator; that is, for $r >0$, define
$$
\omega_r(f;t)_{W_{\kappa,\mu}^T,p} : = \sup_{\theta \le t} 
   \|(T_\theta(W_{\kappa,\mu}^T) - I)^{r/2} f \|_{W_{\kappa,\mu}^T,p}.
$$
Evidently it is related to the modulus of smoothness on the unit ball $B^d$. 
The following relation follows immediately from Proposition \ref{prop:4.4}
and \eqref{eq:4.1}.

\bigskip

\begin{proposition} \label{prop:4.6}
For $f\in L^1(W_{\kappa,\mu}^T)$, 
$$
  \omega_r(f;t)_{W_{\kappa,\mu}^T,p} = 
          \omega_r(f\circ \psi;t)_{W_{\kappa,\mu}^B,p}. 
$$   
\end{proposition}  

Properties of $\omega_r(f;t)_{W_{\kappa,\mu}^T,p}$ can be derived from the 
corresponding ones of $\omega_r(f\circ \psi;t)_{W_{\kappa,\mu}^B,p}$ in 
Proposition \ref{prop:3.6}. We will not write these properties down. The
modulus of smoothness $\omega_r(f;t)_{W_{\kappa,\mu}^T,p}$ is also equivalent 
to the K-functional $K_r(f;t)_{W_{\kappa,\mu}^T,p}$ defined in \cite{X02b}. 
The definition is exactly the same as the one for $W_{\kappa,\mu}^B$,
$$
 K_r(f;t)_{W_{\kappa,\mu}^T,p}:= \inf \left \{ \|f-g\|_{W_{\kappa,\mu}^T,p}+ 
  t^r \|(-D_{\kappa,\mu}^T)^{r/2} g\|_{W_{\kappa,\mu}^T,p} \right \},
$$
where the infimum is taken over all $g \in \CW_r^p(W_{\kappa,\mu}^T)$. The
space $\CW_r^p(W_{\kappa,\mu}^T)$ is defined as its counterpart on $B^d$.

\bigskip

\begin{theorem} \label{thm:4.7}
For $f \in L^p(W_{\kappa,\mu}^T)$, $1 \le p \le \infty$, 
$$
  c_1 \omega_r(f;t)_{W_{\kappa,\mu}^T,p} \le K_r(f;t)_{W_{\kappa,\mu}^T,p} 
    \le c_2 \omega_r(f;t)_{W_{\kappa,\mu}^T,p}.  
$$
\end{theorem}

\begin{proof}
Because of Proposition \ref{prop:4.6} and Proposition \ref{thm:3.7}, it 
suffices to show that 
$$
K_r(f;t)_{W_{\kappa,\mu}^T,p} =  K_r(f\circ \psi;2 t)_{W_{\kappa,\mu}^B,p}.
$$
By \eqref{eq:4.9} and \eqref{eq:4.1}
\begin{align*}
K_r(f;t)_{W_{\kappa,\mu}^T,p} & =
 \inf_g \left \{ \|f\circ \psi - g\circ \psi\|_{W_{\kappa,\mu}^B,p} + 
  2^r t^r \|(-D_{\kappa,\mu}^B)^{r/2} (g\circ \psi)\|_{W_{\kappa,\mu}^B,p}
     \right \}\\
 &  = \inf_{g_0} \left \{ \|f\circ \psi - g_0\|_{W_{\kappa,\mu}^B,p} + 
  2^r t^r \|(-D_{\kappa,\mu}^B)^{r/2} g_0\|_{W_{\kappa,\mu}^B,p} \right \}\\
 & := K_r^*(f\circ \psi;t)_{W_{\kappa,\mu}^B,p},
\end{align*}
where the infimum is taken over all $g_0$ such that $g_0=g \circ \psi \in 
\CW_r^p(W_{\kappa,\mu}^B)$. The definition clearly shows that 
$$
K_r(f;t)_{W_{\kappa,\mu}^T,p} =
 K_r^*(f,t)_{W_{\kappa,\mu}^B,p} \ge K_r(f\circ \psi;2 t)_{W_{\kappa,\mu}^B,p}.
$$
We prove that the reverse inequality holds. For any $\delta > 0$, fix a $g \in
\CW_r^p(W_{\kappa,\mu}^B)$ such that 
$$
 K_r(f;2 t)_{W_{\kappa,\mu}^T,p} \ge \|f-g\|_{W_{\kappa,\mu}^B,p} + 
 2^r t^r \|(-D_{\kappa,\mu}^B)^{r/2} g\|_{W_{\kappa,\mu}^B,p} - \delta. 
$$
Let $g_0 (x) = 2^{-d}\sum_{\varepsilon \in \ZZ_2^d} R(\varepsilon)g(x)$, where 
$R(\varepsilon) g(x):= g(\varepsilon x)$ for $\varepsilon \in \ZZ_2^d$. Then
$g_0$ is even in each of its variables. We claim that 
$R(\varepsilon)D_{\kappa,\mu}^B = D_{\kappa,\mu}^B R(\varepsilon)$.
Indeed, since $h_\kappa$ is even for each of its variables, it is invariant 
under $\ZZ_2^d$, so that $R(\varepsilon)\Delta_h = \Delta_h R(\varepsilon)$ 
and, furthermore, 
$$
  \langle x, \nabla \rangle  R(\varepsilon) g (x)= \sum x_i \varepsilon_i 
\partial_i g(\varepsilon x) = R(\varepsilon)  \langle x, \nabla \rangle g(x),
$$
the claimed equality follows from the definition of $D_{\kappa,\mu}^B$. It 
follows that 
$$
 \|(-D_{\kappa,\mu}^B)^{r/2} g_0\|_{W_{\kappa,\mu}^B,p} \le 2^{-d} 
    \sum \|(-D_{\kappa,\mu}^B)^{r/2} R(\varepsilon)g\|_{W_{\kappa,\mu}^B,p} 
  \le  \|(-D_{\kappa,\mu}^B)^{r/2} g\|_{W_{\kappa,\mu}^B,p}.
$$
Clearly, we also have 
$$
 \|f\circ \psi - g_0\|_{W_{\kappa,\mu}^B,p} \le 2^{-d}\sum
   \|f \circ \psi - R(\varepsilon)g\|_{W_{\kappa,\mu}^B,p} =
    \|f\circ \psi - g\|_{W_{\kappa,\mu}^B,p}.
$$
Consequently, since $g_0$ is even in each of its variables and $g_0 \in 
\CW_r^p(W_{\kappa,\mu}^B)$, it follows that 
\begin{align*}
 K_r^*(f;t)_{W_{\kappa,\mu}^T,p} & \le \|f-g_0\|_{W_{\kappa,\mu}^B,p} + 
  2^r t^r \|(-D_{\kappa,\mu}^B)^{r/2} g_0\|_{W_{\kappa,\mu}^B,p}  \\
& \le  \|f-g\|_{W_{\kappa,\mu}^B,p} + 
  2^r t^r \|(-D_{\kappa,\mu}^B)^{r/2} g\|_{W_{\kappa,\mu}^B,p}  \\
& \le  K_r(f;2 t)_{W_{\kappa,\mu}^T,p} + \delta.
\end{align*}
Since $\delta$ is arbitrary, this completes the proof.
\end{proof}

Again, the above equivalence allows us to state the following important
property of the modulus of smoothness. 

\bigskip

\begin{corollary} \label{cor:4.8}
For $f \in L^p(W_{\kappa,\mu}^T)$, $1 \le p \le \infty$, 
$$
\omega_r(f,\delta t)_{W_{\kappa,\mu}^T,p}  
 \le  c\, \max\{1,\delta^r\} \omega_r(f, t)_{W_{\kappa,\mu}^T,p}.
$$
\end{corollary}

Furthermore, we can state the direct and the inverse theorems for the best 
approximation by polynomials in $L^p(W_{\kappa,\mu}^T)$, given in terms of 
the K-functional in \cite{X02b}, in terms of the modulus of smoothness. 
 
\bigskip

\begin{theorem}
For $f \in L^p(W_{\kappa,\mu}^T)$, $1 \le p \le \infty$, 
$$
 E_n(f)_{W_{\kappa,\mu}^T,p} \le c\, \omega_r(f;n^{-1})_{W_{\kappa,\mu}^T,p}.
$$
On the other hand, 
$$
\omega_r(f;n^{-1})_{W_{\kappa,\mu}^T,p} \le c\, n^{-r}\sum_{k=0}^n (k+1)^{r-1} 
    E_k(f)_{W_{\kappa,\mu}^T,p}.
$$
\end{theorem} 

Let us point out that in the case of $d=1$, 
$\omega_r(f;t)_{W_{\kappa,\mu}^T,p}$ is a modulus of smoothness for the Jacobi
weight function $w_{\kappa,\mu}$ on $[-1,1]$. However, it is different from 
the modulus of smoothness defined in the literature using the generalized
translation operator $S_\theta$ (see, for example, \cite{AW,Bav} for $r=1$). 
Using the convolution structure, we can write 
$$
  \proj_{\kappa,\mu} f = f \star_{\kappa,\mu}^T q_n, \qquad
    q_n(t) = p_n^{(\lambda+\mu-1/2,-1/2)}(1)p_n^{(\lambda+\mu-1/2,-1/2)}(t)
$$
Hence, all summation methods of orthogonal expansions with respect to 
$W_{\kappa,\mu}^B$ can be written in the form of $f\star_{\kappa,\mu}^B g_r$, 
where $g_r$ is the same summation method applies to the Jacobi series with 
$(\alpha,\beta) = (\lambda+\mu-1/2,-1/2)$. Consequently, the modulus of 
smoothness can be used to study the summability of orthogonal expansions on
the simplex.


\begin{thebibliography}{99} 

\bibitem{Rub}
        I. A. Aliev and B. Rubin,
        Spherical harmonics associated to the Laplace-Bessel operator 
        and generalized spherical convolutions,
        \textit{Anal. Appl. (Singap.)} \textbf{1} (2003), 81--109.

\bibitem{Askey}
        R. Askey,
        \textit{Orthogonal polynomials and special functions}, 
	Regional Conference Series in Applied Mathematics \textbf{21}, 
        SIAM, Philadelphia, 1975.

\bibitem{AW}
        R. Askey and S. Wainger,
        A convolution structure for Jacobi series, 
        \textit{Amer. J. Math.} \textbf{91} (1969), p. 463-485.

\bibitem{AF} 
	P. Appell and J. K. de F\'eriet.
	\textit{Fonctions hyperg\'eom\'etriques et hypersph\'eriques, 
        Polynomes d'Hermite}, Gauthier-Villars, Paris, 1926.

\bibitem{Bav} 
        H. Bavinck
        \textit{Jacobi series and approximation},
        Mathematical Centre Tracts, No. 39. 
        Mathematisch Centrum, Amsterdam, 1972.

\bibitem{BBP} 
        H. Berens, P. L. Butzer and S. Pawelke, 
        Limitierungsverfahren von Reihen mehrdimensionaler Kugelfunktionen
        und deren Saturationsverhalten, 
        \textit{Publ. Res. Inst. Math. Sci. Ser. A.} \textbf{4} (1968), 
        201-268.

\bibitem{BL}
        H. Berens and Luoqing Li,         
        On the de la Vall\'ee Poussin means on the sphere,
        \textit{Results in Math.}, \textbf{24} (1993), 12-26.         

\bibitem{Berg}
        C. Berg,
        Corps convexes et potentiels sphériques. (French) 
        \textit{Mat.-Fys. Medd. Danske Vid. Selsk.} \textbf{37} (1969), 64 pp. 
     
\bibitem{BC}
        A. Bonami and J-L. Clerc,
        Sommes de Ces\`aro et multiplicateurs des d\'eveloppe-\allowbreak ments
        en harmoniques sph\'eriques,
        \textit{Trans. Amer. Math. Soc.} \textbf{183} (1973), 223-263.

\bibitem{Bo}
        S. Bochner, 
        Positive zonal functions on spheres, 
        \textit{Proc. Nat. Acad. Sci.,}, \textbf{40} (1954), 1141-1147.

\bibitem{Butz}
        P. L. Butzer,
        Legendre transform methods in the solution of basic problems in 
        algebraic approximation,
        \textit{Functions, series, operators, Vol. I, II (Budapest, 1980)},
         277--301, Colloq. Math. Soc. János Bolyai, \textbf{35}, 
        North-Holland, Amsterdam, 1983. 

\bibitem{CZ} 
        A. P. Calderon and A. Zygmund,
        On a problem of Mihlin,
        \textit{Trans. Amer. Math. Soc.}, \textbf{78} (1955), 209-224. 

\bibitem{D1} 
        C. F. Dunkl,  
	Differential-difference operators associated to reflection groups,
        \textit{Trans. Amer. Math. Soc.} \textbf{311} (1989), 167--183.

\bibitem{DX}
        C. F. Dunkl and Yuan Xu,
        \textit{Orthogonal polynomials of several variables},
        Cambridge Univ. Press, 2001. 
 
\bibitem{Er}
	A. Erd\'elyi, W. Magnus, F. Oberhettinger, and F. G. Tricomi,
	\textit{Higher transcendental functions}, 
	McGraw-Hill, New York, 1953.

\bibitem{G} 
        M. I. Ganzburg, 
        Polynomial approximation on the $m$-dimensional ball, 
        in \textit{Approximation Theory IX, Vol. I},       
        Vanderbilt Univ. Press, Nashville, TN, 1998, p. 141-148.        
 
\bibitem{Ge} 
        I. M. Gelfand, 
        Spherical functions in symmetric Riemann spaces,
        \textit{Dokl. Akad. Nauk. SSSR.} \textbf{70} (1950), 5-8.

\bibitem{Kam}
        A. I. Kamzolov,
        The best approximation on the classes of functions $W_p^\alpha(S^n)$ by
        polynomials in spherical harmonics,
        \textit{Mat. Zametki}, \textbf{32} (1982), 285--293; English transl in 
        \textit{Math Notes}, \textbf{32} (1982), 622-628.

\bibitem{LW}
        Luoqing Li and Kunyang Wang,
        \textit{Harmonic analysis and approximation on the unit sphere}
        Science Press, Beijing, 2000.
 
\bibitem{LX}
        Zh.-K, Li and Yuan Xu,
        Summability of orthogonal expansions of several variables, 
        \- {\it J. Approx. Theory}, \textbf{122} (2003), 267-333. 

\bibitem{LN}
        P. I. Lizorkin and S. M. Nikolskii,
        Approximation theory on the sphere,
        \textit{Proc. Steklov Inst. Math.}, \textbf{172} (1987),       
        295-302.         
 
\bibitem{P}
        S. Pawelke, 
        \"Uber Approximationsordnung bei Kugelfunktionen und {\allowbreak}
        algebraischen
        Polynomen, 
        \textit{T\^ohoku Math. J.}, \textbf{24} (1972), 473-486.

\bibitem{Raf}
        S. Rafalson,
        An extremal relation of the theory of approximation of functions by 
        algebraic polynomials,
        \textit{J. Approx. Theory}, \textbf{110}  (2001), 146--170.

\bibitem{Ra}
        D. L. Ragozin, 
        Constructive polynomial approximation on spheres and projective
        spaces,
        \textit{Trans. Amer. Math. Soc.}, \textbf{162} (1971), 157-170.
     
\bibitem{Ros} 
	M. R\"osler,
	Positivity of Dunkl's intertwining operator,
	\textit{Duke Math. J.\/},  \textbf{98} (1999), 445--463.

\bibitem{Rus} 
        Kh. Rustamov,  
        On approximation of functions on the sphere,  
        \textit{Russian Acad. Sci. Izv. Math.}, \textbf{43} (1994), 311-329.

\bibitem{Rus2} 
        Kh. Rustamov,  
        On the equivalence of different moduli of smoothness on the sphere,
        \textit{Proc. Steklov Inst. Math.} \textbf{204} (1994), 
        no. 3, 235--260.  

\bibitem{SW} 
        R. L. Stens and M. Wehrens, 
        Legendre transform methods and best algebraic approximation,
        \textit{Comment. Math. Prace Mat.} \textbf{21} (1980), 351--380. 

\bibitem{Szego}
	G. Szeg\"{o},
	\textit{Orthogonal Polynomials},  
	Amer. Math. Soc. Colloq. Publ. Vol.23, Providence, 4th edition,
        1975.

\bibitem{T}
        K. Trimèche, 
        The Dunkl intertwining operator on spaces of functions and 
        distributions and integral representation of its dual,
        \textit{Integral Transform. Spec. Funct.}  \textbf{12} (2001),
        349--374. 
 
\bibitem{W} 
        M. Wherens,  
        Best approximation on the unit sphere in $\RR^k$, 
        in \textit{Functional Analysis and Approximation} (Oberwolfach, 1980), 
        Birkhauser, Basel, 1981, 233-245. 

\bibitem{X97b}
	Yuan Xu,
	Integration of the intertwining operator for $h$-harmonic polynomials
	associated to reflection groups, 
	\textit{Proc. Amer. Math. Soc.} \textbf{125} (1997), 2963--2973.




\bibitem{X00}
	Yuan Xu,            
	Funk-Hecke formula for orthogonal polynomials on spheres and on balls,
        \textit{Bull. London Math. Soc.} \textit{32} (2000), 447-457.  

\bibitem{X01a}
	Yuan Xu,            
	Orthogonal polynomials and summability in Fourier orthogonal 
        series on spheres and on balls,  
        \textit{Math. Proc. Cambridge Phil. Soc.}, \textbf{31} (2001), 
        139-155.    

\bibitem{X01b} 
	Yuan Xu,            
	Orthogonal polynomials on the ball and on the simplex for weight
        functions with reflection symmetries,
        \textit{Constr. Approx.}, \textbf{17} (2001), 383-412.

\bibitem{X02a} 
	Yuan Xu,            
        Approximation by means of $h$-harmonic polynomials on the unit sphere,
        \textit{Adv. in Comp. Math}, to appear.\\
        \texttt{http://math.uoregon.edu/\~{}yuan} (ApproxSph.ps.gz). 

\bibitem{X02b} 
	Yuan Xu,            
        Weighted approximation of functions on the unit sphere,
        \textit{Const. Approx.}, to appear.
        \texttt{http://math.uoregon.edu/\~{}yuan} (BestApp.ps.gz). 

\bibitem{X03} 
	Yuan Xu,            
        Almost everywhere convergence of orthogonal expansions of several 
        variables, 
        \textit{Const. Approx.}, to appear.         
        \texttt{http://math.uoregon.edu/\~{}yuan} (Almost.ps.gz). 

\end{thebibliography}
\enddocument